\documentclass[11pt, a4paper]{article}
\usepackage[cp1251]{inputenc}
\usepackage{amsmath} \usepackage{euscript} \usepackage{marvosym}
\usepackage[dvipsnames]{xcolor}
\usepackage{mymatrx6}
\usepackage{graphicx}
\usepackage{longtable}
\usepackage{mathrsfs}
\usepackage{amssymb, amscd}
\usepackage{comment}
\def\mymatrix{\MyMatrixwithdelims..}
\Linestrue\Autonumfalse
\begin{document}

\newcounter{bnomer} \newcounter{snomer}
\newcounter{bsnomer}
\setcounter{bnomer}{0}
\renewcommand{\thesnomer}{\thebnomer.\arabic{snomer}}
\renewcommand{\thebsnomer}{\thebnomer.\arabic{bsnomer}}
\renewcommand{\refname}{\begin{center}\large{\textbf{References}}\end{center}}

\setcounter{MaxMatrixCols}{14}

\newcommand\restr[2]{{
  \left.\kern-\nulldelimiterspace 
  #1 
  \right|_{#2} 
}}

\newcommand{\sect}[1]{%
\setcounter{snomer}{0}\setcounter{bsnomer}{0}
\refstepcounter{bnomer}
\par\bigskip\begin{center}\large{\textbf{\arabic{bnomer}. {#1}}}\end{center}}
\newcommand{\sst}[1]{%
\refstepcounter{bsnomer}
\par\bigskip\textbf{\arabic{bnomer}.\arabic{bsnomer}. {#1}}\par}
\newcommand{\defi}[1]{%
\refstepcounter{snomer}
\par\medskip\textbf{Definition \arabic{bnomer}.\arabic{snomer}. }{#1}\par\medskip}
\newcommand{\theo}[2]{%
\refstepcounter{snomer}
\par\textbf{Theorem \arabic{bnomer}.\arabic{snomer}. }{#2} {\emph{#1}}\hspace{\fill}$\square$\par}
\newcommand{\mtheop}[2]{%
\refstepcounter{snomer}
\par\textbf{Theorem \arabic{bnomer}.\arabic{snomer}. }{\emph{#1}}
\par\textsc{Proof}. {#2}\hspace{\fill}$\square$\par}
\newcommand{\mcorop}[2]{%
\refstepcounter{snomer}
\par\textbf{Corollary \arabic{bnomer}.\arabic{snomer}. }{\emph{#1}}
\par\textsc{Proof}. {#2}\hspace{\fill}$\square$\par}
\newcommand{\mtheo}[1]{%
\refstepcounter{snomer}
\par\medskip\textbf{Theorem \arabic{bnomer}.\arabic{snomer}. }{\emph{#1}}\par\medskip}
\newcommand{\theobn}[1]{%
\par\medskip\textbf{Theorem. }{\emph{#1}}\par\medskip}
\newcommand{\theoc}[2]{%
\refstepcounter{snomer}
\par\medskip\textbf{Theorem \arabic{bnomer}.\arabic{snomer}. }{#1} {\emph{#2}}\par\medskip}
\newcommand{\mlemm}[1]{%
\refstepcounter{snomer}
\par\medskip\textbf{Lemma \arabic{bnomer}.\arabic{snomer}. }{\emph{#1}}\par\medskip}
\newcommand{\mprop}[1]{%
\refstepcounter{snomer}
\par\medskip\textbf{Proposition \arabic{bnomer}.\arabic{snomer}. }{\emph{#1}}\par\medskip}
\newcommand{\theobp}[2]{%
\refstepcounter{snomer}
\par\textbf{Theorem \arabic{bnomer}.\arabic{snomer}. }{#2} {\emph{#1}}\par}
\newcommand{\theop}[2]{%
\refstepcounter{snomer}
\par\textbf{Theorem \arabic{bnomer}.\arabic{snomer}. }{\emph{#1}}
\par\textsc{Proof}. {#2}\hspace{\fill}$\square$\par}
\newcommand{\theosp}[2]{%
\refstepcounter{snomer}
\par\textbf{Theorem \arabic{bnomer}.\arabic{snomer}. }{\emph{#1}}
\par\textsc{Sketch of the proof}. {#2}\hspace{\fill}$\square$\par}
\newcommand{\exam}[1]{%
\refstepcounter{snomer}
\par\medskip\textbf{Example \arabic{bnomer}.\arabic{snomer}. }{#1}\par\medskip}
\newcommand{\deno}[1]{%
\refstepcounter{snomer}
\par\textbf{Notation \arabic{bnomer}.\arabic{snomer}. }{#1}\par}
\newcommand{\lemm}[1]{%
\refstepcounter{snomer}
\par\textbf{Lemma \arabic{bnomer}.\arabic{snomer}. }{\emph{#1}}\hspace{\fill}$\square$\par}
\newcommand{\lemmp}[2]{%
\refstepcounter{snomer}
\par\medskip\textbf{Lemma \arabic{bnomer}.\arabic{snomer}. }{\emph{#1}}
\par\textsc{Proof}. {#2}\hspace{\fill}$\square$\par\medskip}
\newcommand{\coro}[1]{%
\refstepcounter{snomer}
\par\textbf{Corollary \arabic{bnomer}.\arabic{snomer}. }{#1}\hspace{\fill}$\square$\par}
\newcommand{\mcoro}[1]{%
\refstepcounter{snomer}
\par\textbf{Corollary \arabic{bnomer}.\arabic{snomer}. }{\emph{#1}}\par\medskip}
\newcommand{\corop}[2]{%
\refstepcounter{snomer}
\par\textbf{Corollary \arabic{bnomer}.\arabic{snomer}. }\emph{#1}
\par\textsc{Proof}. {#2}\hspace{\fill}$\square$\par}
\newcommand{\nota}[1]{%
\refstepcounter{snomer}
\par\medskip\textbf{Remark \arabic{bnomer}.\arabic{snomer}. }{#1}\par\medskip}
\newcommand{\propp}[2]{%
\refstepcounter{snomer}
\par\medskip\textbf{Proposition \arabic{bnomer}.\arabic{snomer}. }{\emph{#1}}
\par\textsc{Proof}. {#2}\hspace{\fill}$\square$\par\medskip}
\newcommand{\hypo}[1]{%
\refstepcounter{snomer}
\par\medskip\textbf{Conjecture \arabic{bnomer}.\arabic{snomer}. }{\emph{#1}}\par\medskip}
\newcommand{\prop}[1]{%
\refstepcounter{snomer}
\par\textbf{Proposition \arabic{bnomer}.\arabic{snomer}. }{\emph{#1}}\hspace{\fill}$\square$\par}

\newcommand{\proof}[2]{%
\par\medskip\textsc{Proof{#1}}. \hspace{-0.2cm}{#2}\hspace{\fill}$\square$\par\medskip}

\makeatletter
\def\iddots{\mathinner{\mkern1mu\raise\p@
\vbox{\kern7\p@\hbox{.}}\mkern2mu
\raise4\p@\hbox{.}\mkern2mu\raise7\p@\hbox{.}\mkern1mu}}
\makeatother

\newcommand{\okr}[2]{%
\refstepcounter{snomer}
\par\medskip\textbf{{#1} \arabic{bnomer}.\arabic{snomer}. }{\emph{#2}}\par\medskip}

\newcommand{\Ind}[3]{%
\mathrm{Ind}_{#1}^{#2}{#3}}
\newcommand{\Res}[3]{%
\mathrm{Res}_{#1}^{#2}{#3}}
\newcommand{\epsi}{\varepsilon}
\newcommand{\tri}{\triangleleft}
\newcommand{\Supp}[1]{%
\mathrm{Supp}(#1)}
\newcommand{\SSu}[1]{%
\mathrm{SingSupp}(#1)}

\newcommand{\gee}{\geqslant}
\newcommand{\reg}{\mathrm{reg}}
\newcommand{\Dyn}{\mathrm{Dyn}}
\newcommand{\Ann}{\mathrm{Ann}\,}
\newcommand{\Cent}[1]{\mathbin\mathrm{Cent}({#1})}
\newcommand{\PCent}[1]{\mathbin\mathrm{PCent}({#1})}
\newcommand{\Irr}[1]{\mathbin\mathrm{Irr}({#1})}
\newcommand{\Exp}[1]{\mathbin\mathrm{Exp}({#1})}
\newcommand{\empr}[2]{[-{#1},{#1}]\times[-{#2},{#2}]}
\newcommand{\sreg}{\mathrm{sreg}}
\newcommand{\ilm}{\varinjlim}
\newcommand{\wdth}{\mathrm{wd}}
\newcommand{\plm}{\varprojlim}
\newcommand{\codim}{\mathrm{codim}\,}
\newcommand{\GKdim}{\mathrm{GKdim}\,}
\newcommand{\chara}{\mathrm{char}\,}
\newcommand{\rk}{\mathrm{rk}\,}
\newcommand{\chr}{\mathrm{ch}\,}
\newcommand{\Ker}{\mathrm{Ker}\,}
\newcommand{\id}{\mathrm{id}}
\newcommand{\Ad}{\mathrm{Ad}}
\newcommand{\Gh}{\mathrm{Gh}}
\newcommand{\col}{\mathrm{col}}
\newcommand{\row}{\mathrm{row}}
\newcommand{\high}{\mathrm{high}}
\newcommand{\low}{\mathrm{low}}
\newcommand{\pho}{\hphantom{\quad}\vphantom{\mid}}
\newcommand{\fho}[1]{\vphantom{\mid}\setbox0\hbox{00}\hbox to \wd0{\hss\ensuremath{#1}\hss}}
\newcommand{\wt}{\widetilde}
\newcommand{\wh}{\widehat}
\newcommand{\ad}[1]{\mathrm{ad}_{#1}}
\newcommand{\tr}{\mathrm{tr}\,}
\newcommand{\GL}{\mathrm{GL}}
\newcommand{\SL}{\mathrm{SL}}
\newcommand{\SO}{\mathrm{SO}}
\newcommand{\Or}{\mathrm{O}}
\newcommand{\Sp}{\mathrm{Sp}}
\newcommand{\Sa}{\mathrm{S}}
\newcommand{\Ua}{\mathrm{U}}
\newcommand{\Andre}{\mathrm{Andre}}
\newcommand{\Aord}{\mathrm{Aord}}
\newcommand{\Mat}{\mathrm{Mat}}
\newcommand{\Stab}{\mathrm{Stab}}
\newcommand{\htt}{\mathfrak{h}}
\newcommand{\spt}{\mathfrak{sp}}
\newcommand{\slt}{\mathfrak{sl}}
\newcommand{\sot}{\mathfrak{so}}

\newcommand{\vfi}{\varphi}
\newcommand{\aad}{\mathrm{ad}}
\newcommand{\vpi}{\varpi}
\newcommand{\teta}{\vartheta}
\newcommand{\Bfi}{\Phi}
\newcommand{\Fp}{\mathbb{F}}
\newcommand{\Rp}{\mathbb{R}}
\newcommand{\Zp}{\mathbb{Z}}
\newcommand{\Cp}{\mathbb{C}}
\newcommand{\Ap}{\mathbb{A}}
\newcommand{\Pp}{\mathbb{P}}
\newcommand{\Kp}{\mathbb{K}}
\newcommand{\Np}{\mathbb{N}}
\newcommand{\ut}{\mathfrak{u}}
\newcommand{\at}{\mathfrak{a}}
\newcommand{\glt}{\mathfrak{gl}}
\newcommand{\hei}{\mathfrak{hei}}
\newcommand{\nt}{\mathfrak{n}}
\newcommand{\kt}{\mathfrak{k}}
\newcommand{\mt}{\mathfrak{m}}
\newcommand{\rt}{\mathfrak{r}}
\newcommand{\rad}{\mathfrak{rad}}
\newcommand{\bt}{\mathfrak{b}}
\newcommand{\unt}{\underline{\mathfrak{n}}}
\newcommand{\gt}{\mathfrak{g}}
\newcommand{\vt}{\mathfrak{v}}
\newcommand{\pt}{\mathfrak{p}}
\newcommand{\Xt}{\mathfrak{X}}
\newcommand{\Po}{\mathcal{P}}
\newcommand{\PV}{\mathcal{PV}}
\newcommand{\Uo}{\EuScript{U}}
\newcommand{\Fo}{\EuScript{F}}
\newcommand{\Do}{\EuScript{D}}
\newcommand{\Eo}{\EuScript{E}}
\newcommand{\Jo}{\EuScript{J}}
\newcommand{\Iu}{\mathcal{I}}
\newcommand{\Mo}{\mathcal{M}}
\newcommand{\Nu}{\mathcal{N}}
\newcommand{\Ro}{\mathcal{R}}
\newcommand{\Co}{\mathcal{C}}
\newcommand{\Ko}{\mathcal{K}}
\newcommand{\So}{\mathcal{S}}
\newcommand{\Lo}{\mathcal{L}}
\newcommand{\Ou}{\mathcal{O}}
\newcommand{\Uu}{\mathcal{U}}
\newcommand{\Tu}{\mathcal{T}}
\newcommand{\Au}{\mathcal{A}}
\newcommand{\Vu}{\mathcal{V}}
\newcommand{\Du}{\mathcal{D}}
\newcommand{\Bu}{\mathcal{B}}
\newcommand{\Sy}{\mathcal{Z}}
\newcommand{\Sb}{\mathcal{F}}
\newcommand{\Gr}{\mathcal{G}}
\newcommand{\Xu}{\mathcal{X}}
\newcommand{\Op}{\mathbb{O}}
\newcommand{\chv}{\mathrm{chv}}
\newcommand{\rtc}[1]{C_{#1}^{\mathrm{red}}}
\newcommand{\Dsr}{D_\text{sreg}}
\newcommand{\Dr}{D_\text{reg}}
\newcommand{\e}{\epsi}

\author{Mikhail Ignatev\and Mikhail Venchakov}
\date{}
\title{Orbits of submaximal dimension for Sylow $p$-subgroups\\ of finite classical orthogonal groups}\maketitle
\begin{abstract} Let $U$ be a Sylow $p$-subgroup in a classical orthogonal group over a finite field with $q$ elements of characteristic $p$ large enough. The coadjoint orbits of the group $U$ play the key role in the description of irreducible complex characters of $U$. In the paper, we provide a classification of such orbits of submaximal dimension. As a corollary, we compute the number of all orbits mentioned above. It turned out that each of these numbers is a polynomial in $q - 1$ with integer non-negative coefficients, which agrees with Isaacs' conjecture.


\end{abstract}

\sect{Introduction}

\let\thefootnote\relax\footnote{The research was supported by RSF (project No. 25--21--00219), \texttt{https://rscf.ru/en/project/25-21-00219/}.}

Let $U$ be a unipotent algebraic group over a finite field $\Fp_q$ of sufficiently large characteristic $p$. The main tool in representation theory of $U$ is the orbit method created in 1962 by A.A. Kirillov, see~\cite{Kirillov62}, \cite{Kirillov04}, \cite{Kazhdan77}. The key idea of the orbit method says that the irreducible representations of $U$ are in one-to-one correspondence with the coadjoint orbits of this group. Namely, the group $U$ acts on its Lie algebra $\ut$ via the adjoint action: $$g\in U,~x\in\ut\mapsto gxg^{-1}\in\ut;$$ the dual action of $U$ on the dual space $\ut^*$ is called coadjoint. It turns out that there is a natural bijection between the set $\Irr{U}$ of all irreducible complex characters of $U$ and the set $\ut^*/U$ of coadjoint orbits.

Let $U$ be a maximal unipotent subgroup (or, equivalently, a Sylow $p$-subgroup) in a simple classical group over $\Fp_q$. A complete description of all coadjoint orbits for the group $U$ is a wild problem, so a natural question is how to describe certain important classes of orbits and the corresponding irreducible characters. For the case $A_{n-1}$, orbits of maximal possible dimension (we call such orbits regular) were classified in Kirillov's first work on the orbit method \cite{Kirillov62}. All of them are associated with so-called orthogonal rook placements.

For other classical root systems, so-called orbits associated with the Kostant cascades have maximal possible dimension, see \cite{Kostant12} and \cite{Kostant13}, but not all the orbits of maximal dimension have such a form. For $A_{n-1}$, orbits of submaximal dimension (we call such orbits subregular) were classified by A.N. Panov \cite{IgnatevPanov09}. Such orbits also correspond to orthogonal rook placements (modulo adding simple root covectors to the canonical forms on them). A classification of orbits of maximal dimension in type $C_n$ follows from C. Andre and A. Neto's paper \cite{AndreNeto06}, see also \cite{Venchakov}. A description of orbits of maximal possible dimension for the root systems $B_n$ and $D_n$, as well as of orbits of submaximal dimension for $C_n$ was obtained in \cite{IgnatevVenchakov26}. 

\newpage
The main goal of this paper is to classify subregular coadjoint orbits for the root systems $B_n$ and $D_n$. 
The structure of the paper is as follows.
In the next section we introduce the necessary notation and definitions. The main results are formulated in Section \ref{sect:sea_battle_pol}. More precisely, we provide an explicit classification of the subregular orbits (i.e., the orbits of submaximal dimension), prove Theorem~\ref{theo_orb_submax_dim} and calculate their number as a corollary. The main technical tool which was used is the method of C-patterns and C-quatterns, invented in \cite{GoodwinMoschRoehrle16} by S.M. Goodwin, P. Mosch and G. R\"ohrle, and then applied to description of orbits of certain special classes in \cite{IgnatevPetukhov25} by the first author and A. Petukhov, and in \cite{Venchakov} by the second author.

A longstanding G. Higman's conjecture \cite{Higman60} states that the number of conjugacy classes for~$U=U_n(q)$, the group of all upper-triangular matrices over $\Fp_q$ with 1's on the diagonal (which is a Sylow $p$-subgroup in $\SL_n(\Fp_q)$), is a polynomial in~$q$. It is easy to check that the number of conjugacy classes for $U$ coincides with the number of coadjoint $U$-orbits, or, equivalently, with the number of irreducible characters of $U$. Of course, it is interesting to consider an analogue of Higman's conjecture for an arbitrary $U$, not only for $U_n(q)$; we will denote the number of irreducible characters of the group $U$ by $O(q)$. Fourteen years after Higman, G. Lehrer conjectured \cite{Lehrer74} that, for $U=U_n(q)$, even the number $O_e(q)$ of characters of degree $q^e$ is polynomial in $q$. In 2007, I.M. Isaacs \cite{Isaacs07} made a stronger conjecture that, for $U=U_n(q)$, $O_e(q)$ is in fact polynomial in $q-1$ with nonnegative integer coefficients. In the past twenty-five years, significant progress in studying these conjectures has been made for $U=U_n(q)$, see, e.g., \cite{VeraLopezArregi03}, \cite{Marjoram97'}, \cite{Marjoram97}, \cite{IgnatevPanov09}, \cite{Marjoram99}, \cite{Loukaki11}, \cite{Le10}, \cite{Marberg11}.

For Sylow subgroups $U$ in other Chevalley groups, the situation is as follows. In 1999, Marjoram computed $O_{\mu(\Phi)}(q)$ for orthogonal group (i.e., for $\Phi=B_n$ or $D_n$), where $q^{\mu(\Phi)}$ is the maximal possible degree of an irreducible character of $U$, see \cite{Marjoram99}. For the symplectic case (i.e., for $\Phi=C_n$), the formula for $O_{\mu(\Phi)}(q)$ follows from the results of C.A.M. Andre and A.-M. Neto on so-called supercharacters published in 2006 \cite{AndreNeto06}. In 2016, S.M. Goodwin, P. Mosch and G. R\"ohrle calculated $O_e(q)$ and proved Isaacs' conjecture for all possible $e$ and all finite Chevalley groups $G(q)$ of rank $\leq 8$, except $E_8$. Our classification of subregular orbits for $B_n$ and $D_n$ gives a proof of Isaacs' conjecture in that case.

We thank Alexey Petukhov for very useful discussions.


\sect{Main definitions} \label{main_definitions}
To begin with we present some basic facts about root systems of simple algebras. We will denote root systems of type $B_n$ or $D_n$ by $\Phi$. As usual, we will identify them with the following subsets of $\Rp^n$:
$$
\begin{array}{ll}
&B_n = \{\pm\epsi_i\pm\epsi_j, 1\leq i < j\leq
n\}\cup\{\pm\epsi_i, 1\leq i\leq n\},\\
&D_n = \{\pm\epsi_i\pm\epsi_j, 1\leq i < j\leq n\}.\\
\end{array}
$$
where $\{\epsi_i\}_{i=1}^n$ is the standard basis in $\Rp^n$. Pick the set of simple roots $\Delta=\Delta(\Phi)$ as in
\cite{Bourbaki03}:
$$
\begin{array}{ll}
&\Delta(B_n)=\{\epsi_i-\epsi_{i+1},1\leq i\leq n-1\}\cup\{\epsi_n\},\\
&\Delta(D_n)=\{\epsi_i-\epsi_{i+1},1\leq i\leq n-1\}\cup\{\epsi_{n-1}+\epsi_n\}.\\
\end{array}
$$
The set of positive roots $\Phi^+\supset\Delta$ looks as follows:
$$
\begin{array}{ll}
&B_n^+ = \{\epsi_i\pm\epsi_j, 1\leq i < j\leq n\}\cup\{\epsi_i,
1\leq i\leq n\},\\
&D_n^+ = \{\epsi_i\pm\epsi_j, 1\leq i < j\leq n\}.\\
\end{array}
$$
Let
$$
m=\left\{\begin{array}{ll}2n+1,&\mbox{if }\Phi=B_n,\\
2n,&\mbox{if }\Phi=D_n.
\end{array}\right.
$$
Denote by $\ut=\ut(\Phi)$ the subalgebra of $\mathfrak{gl}_m(\Fp_q)$
spanned by the vectors $e_{\alpha}$, $\alpha\in\Phi^+$, where
$$
\begin{array}{ll}
&e_{\epsi_i}=e_{0,i}-e_{-i,0},\quad 1\leq i\leq n,\\
&e_{\epsi_i-\epsi_j}=e_{j,i}-e_{-i,-j},\quad 1\leq i<j\leq n,\\
&e_{\epsi_i+\epsi_j}=e_{-j,i}-e_{-i,j},
\quad 1\leq i<j\leq n.
\end{array}
$$
Here we number the rows and the columns of $m\times m$ matrix by the indices $$1,2,\ldots,n,0,-n,\ldots,-2,-1$$ (there is no index $0$ in the case $D_n$), and we denote by $e_{a, b}$ the usual elementary matrix. Of course, it is a maximal nilpotent subalgebra in the corresponding classical algebra $\gt= \gt(\Phi)$. In particular, $\dim\ut = |\Phi^+|$. \exam{Here we schematically drew the algebras $B_3$ and $D_4$. The convention is as follows: given $1\leq j<i$, the square $(i,j)$ (respectively, $(-i,j)$, $(0,j)$) corresponds to the root $\epsi_j-\epsi_i$ (respectively, $\epsi_j+\epsi_i$, $\epsi_j$).
\begin{center}\small
$\mymatrix{
\Note{1}\lNote{1}\pho& \Note{2}\pho& \Note{3}\pho& \Note{0}\pho& \Note{-3}\pho& \Note{-2}\pho & \Note{-1}\pho\\
\lNote{2}\Top{2pt}\Rt{2pt} \pho & \pho& \pho& \pho& \pho& \pho & \pho\\
\lNote{3}\pho & \Top{2pt}\Rt{2pt} \pho& \pho& \pho& \pho& \pho & \pho\\
\lNote{0}\pho & \pho& \Top{2pt}\Rt{2pt} \pho& \pho& \pho& \pho & \pho\\
\lNote{-3}\pho & \pho& \Top{2pt}\Lft{2pt} & \Top{2pt}\Rt{2pt} \pho& \pho& \pho & \pho\\
\lNote{-2}\pho& \Top{2pt}\Lft{2pt} & \pho& \pho& \Top{2pt}\Rt{2pt} \pho& \pho & \pho\\
\lNote{-1}\Top{2pt} & \pho& \pho& \pho& \pho& \Top{2pt}\Rt{2pt}\pho & \pho\\
}\quad\mymatrix{
\lNote{1}\Note{1}\pho& \Note{2}\pho& \Note{3}\pho& \Note{4}\pho& \Note{-4}\pho& \Note{-3}\pho & \Note{-2}\pho & \Note{-1}\pho\\
\lNote{2}\Top{2pt}\Rt{2pt} \pho & \pho& \pho& \pho& \pho& \pho & \pho & \pho\\
\lNote{3}\pho & \Top{2pt}\Rt{2pt} \pho& \pho& \pho& \pho& \pho & \pho & \pho\\
\lNote{4}\pho & \pho& \Top{2pt}\Rt{2pt} \pho& \pho& \pho& \pho & \pho & \pho\\
\lNote{-4}\pho &\pho & \pho& \Top{2pt}\Lft{2pt}\Rt{2pt} & \pho& \pho& \pho & \pho\\
\lNote{-3}\pho &\pho& \Top{2pt}\Lft{2pt} & \pho& \Top{2pt}\Rt{2pt}\pho& \pho& \pho & \pho\\
\lNote{-2}\pho &\Top{2pt}\Lft{2pt} & \pho& \pho& \pho& \Top{2pt}\Rt{2pt} \pho& \pho & \pho\\
\lNote{-1}\Top{2pt} & \pho& \pho& \pho& \pho& \pho & \Top{2pt}\Rt{2pt}\pho &\pho\\
}
$
\end{center}
}

Furthermore, we define the functions 
$$
\begin{array}{ll}
&\col\colon\Phi^+\to\{1,\ldots,n\}\colon\col(\epsi_i\pm\epsi_j)=\col(\epsi_i)=i,\\
&\row\colon\Phi^+\to\{-n,\ldots,n\}\colon\row(\epsi_i\pm\epsi_j)=\mp j,\row(\epsi_i)=0.\\
\end{array}
$$
For arbitrary $-n+1\leq i\leq n-1$ and $1\leq j\leq n$, the sets
$$
\begin{array}{ll}
&R_i = R_i(\Phi) = \{\alpha\in\Phi^+\mid \row(\alpha)=i\},\\
&C_j = C_j(\Phi) = \{\alpha\in\Phi^+\mid \col(\alpha)=j\}\\
\end{array}
$$
are called the $i$th\emph{ row} and the $j$th\emph{ column}
$\Phi^+$ respectively. We introduce the mirror order on the 
set of indices $$1\prec2\prec\ldots\prec
n\prec0\prec-n\prec\ldots\prec-2\prec-1,$$ and the following total orders on $\Phi^+$:
\begin{equation*}
\alpha\prec\beta~\stackrel{\mathrm{def}}{\iff}~
\col(\beta)\prec\col(\alpha)~\mbox{or}~
\col(\beta)=\col(\alpha),~\row(\beta)\succ\row(\alpha).
\end{equation*}
For example, for $\Phi=B_6$ we have
$\epsi_2+\epsi_5\succ\epsi_2\succ\epsi_2-\epsi_4\succ\epsi_3-\epsi_6$.

For given $\Phi$ denote the maximal unipotent subgroup $\exp(\ut(\Phi))$ in the corresponding classical finite group $G$ by~$U(\Phi)$ (we will write just $\ut$ and $U$ in cases when $\Phi$ is fixed). In the sequel, we will assume everywhere that $p=\chara{\Fp_q}\gee m$. Under this assumption, the map
$$
\exp\colon\ut\to U\colon x\mapsto \sum_{i=0}^{m-1}\frac{x^i}{i!}
$$
is well-defined and is in fact a bijection (and also an isomorphism of algebraic varieties over $\overline{\Fp_q}$);
we denote the inverse map by $\ln$. Furthermore, the 
Baker--Campbell--Hausdorff formula claims that, for a Lie subalgebra $\at\subset\ut$ and arbitrary $u, v\in\at$, one has
\begin{equation*}
\exp(u)\exp(v)=\exp(u + v + \tau(u,
v)),
\end{equation*} where
$\tau(u, v)\in[\at, \at]$ (here $[\at, \at]=\langle[x, y], x,
y\in\at\rangle_{\Fp_q}$).

The group $U$ acts on its Lie algebra $\ut$ by the adjoint action; the dual action of $U$ on the $\Fp_q$-dual space $\ut^*$ is called \emph{coadjoint}. Using the non-degenerate form $\langle A, B\rangle=\mathrm{tr}(AB)$ on $\mathfrak{gl}_m(\Fp_q)$, one can identify the dual space $\ut^*$ with the space $\ut^t$ (in this case, $e^*_{\alpha}=e_{\alpha}^T/2$). Under this identification, the coadjoint action has the following form:
$$
g.x = \mathrm{pr}(gxg^{-1}),~g\in U,~x\in\ut^*
$$
(here we denote the projection
$\mathfrak{gl}_m(\Fp_q)\to\ut^*$ by $\mathrm{pr}$).
\defi{Let $D=\{\beta_1,\ldots,\beta_t\}\subset\Phi^+$ be an arbitrary subset of the positive roots, and
$\xi=(\xi_{\beta})_{\beta\in D}$ be a set of non-zero constants from $\Fp_q$. Put
$$
f=f_{D, \xi} = \sum_{\beta\in D}\xi_{\beta}e_{\beta}^*\in\ut^*. $$ We say that the coadjoint orbit $\Omega=\Omega_{D,
\xi}\subset\ut^*$ of the linear form $f$ is \emph{associated} with $D$, and $f$ is called the
\emph{canonical form} on $\Omega$.}

\defi{Let $\beta\in\Phi^+$. Roots $\alpha$, $\gamma\in\Phi^+$ are called
$\beta$-\emph{singular} if $\alpha+\gamma=\beta$. The set of all $\beta$-singular roots is denoted by $S(\beta)$ (see
\cite{Andre95ii} for $A_n$ and \cite{AndreNeto06} for
$B_n$, $D_n$).} It is easy to see that the singular roots look as follows:
\begin{equation*}
\begin{split}
&S(\epsi_i-\epsi_j)=\bigcup_{l=i+1}^{j-1}\{\epsi_i-\epsi_l,\epsi_l-\epsi_j\},~1\leq i < j\leq n,\\
&S(\epsi_i) = \bigcup_{l=i+1}^{n}\{\epsi_i-\epsi_l,\epsi_l\},~1\leq i\leq n,\\
&S(\epsi_i+\epsi_j)=\bigcup_{l=i+1}^{j-1}\{\epsi_i-\epsi_l,\epsi_l+\epsi_j\}
\cup\bigcup_{l=j+1}^n\{\epsi_i-\epsi_l,\epsi_j+\epsi_l\}\cup\\
&\bigcup_{l=j+1}^n\{\epsi_i+\epsi_l,\epsi_j-\epsi_l\}\cup S_{ij}
,~1\leq i<j\leq
n,\mbox{ where}\\
&S_{ij}=
\begin{cases}\{\epsi_i, \epsi_j\},&\mbox{if }\Phi=B_n,\\
\varnothing,&\mbox{if }\Phi=D_n.
\end{cases}.
\end{split}\label{formula_sing_roots}
\end{equation*}

\nota{Note that, given $\beta\in\Phi^+$, 
the column of $\beta$ contains exactly one of the roots from each pair of $\beta$-singular roots whose sum is $\beta$. Precisely, put
\begin{equation*}
S^+(\beta)=S(\beta)\cap C_{\text{col}(\beta)}
\end{equation*}
so that $S(\beta)=S^+(\beta)\cup S^-(\beta)$\label{nota_one_row_col}, where $S^-(\beta)=S(\beta)\setminus S^+(\beta)$.}

Next, we need to recall some necessary notation from \cite{IgnatevPetukhov25} which will be used in proofs of the main results of the paper.

Let $X$ be a subset of $\Phi^+$. Set $\ut_X:=\bigoplus_{\alpha\in X}\Fp e_\alpha$ with the (Lie) bracket given by the formula $$[e_\alpha, e_\beta]=\begin{cases}[e_\alpha, e_\beta],& {\rm~if~}\alpha+\beta\in X,\\0,&\text{if }\alpha+\beta\notin X\end{cases}$$ for all $\alpha, \beta\in X$
(note that if $\alpha, \beta, \alpha+\beta\in \Phi$, then $[e_\alpha, e_\beta]$ is proportional to $e_{\alpha+\beta})$. Such a bilinear bracket $[\cdot, \cdot ]$ does not always define a Lie algebra but it does define a Lie algebra under the assumption that $X$ is a $C$-pattern (a shorthand for ``combinatorial counterpart of a pattern subgroup'') or $C$-quattern (a shorthand for ``combinatorial counterpart of a quattern subgroup''), see the next definitions. 

\defi{We will say that $X$ is a {\it C-pattern}  if $X$ satisfies the following condition: $$\text{if }\alpha, \beta\in X\text{ and }\alpha+\beta\in\Phi^+\text{ then }\alpha+\beta\in X.$$
This condition is clearly equivalent to the statement that $\ut_X$ is a sub-Lie algebra of $\ut$.}
\defi{ We will say that $X$ is a {\it C-quattern} if $X=X_+\setminus X_-$ for two C-patterns $X_+, X_-$ such that $X_-\subset X_+$ and
$$\text{if }\alpha_-\in X_-,\alpha_+\in X_+~\text{ and }\alpha_-+\alpha_+\in X_+ \text{ then }\alpha_-+\alpha_+\in X_-.$$
This is clearly equivalent to the statement that $\ut_{X_-}$ is an ideal of the Lie algebra $\ut_{X_+}$ and in this case we have $\ut_X\cong \ut_{X_+}/\ut_{X_-}$ where the right hand side is a quotient of Lie algebra by its ideal and thus the left hand side is a Lie algebra.}
From now on we assume that $X$ is always a C-quattern. We denote by $U_X$ the unipotent group with Lie algebra $\ut_X$. Every C-pattern $X$ is a C-quattern for $X_+=X, X_-={\varnothing}$ and for every C-quattern $X$ the data $(\ut_X, [\cdot, \cdot])$ provides a structure of Lie algebra on $\ut_X$. 

Set 
$$Z(X):=\{\alpha\in X\mid (\alpha+X)\cap X={\varnothing}\}.$$ 
Then it is easy to verify that $\ut_{Z(X)}$ coincides with the center of $\ut_X$. In particular, this implies that $Z(X)\ne\varnothing$ for any non-empty C-quattern $X$.

\defi{Pick $f\in \ut_X^*$ and consider $Z\subset Z(X)$. We will say that $f$ is {\it $Z$-saturated} if $f(e_\alpha)\ne0$ for all $\alpha\in Z$. We will say that $f$ is {\it saturated} if $f$ is $Z(X)$-saturated. 
We denote the variety of $Z$-saturated linear forms by $\ut^*_{X; Z}$.}
It is clear that each $\ut_{X; Z}^*$ is $U_X$-stable and thus is a collection of coadjoint orbits. We map $\ut^*_{X; Z}$ to $\ut^*$ by the formula $e_\alpha^*\mapsto e_\alpha^*$ and denote by $\underline{\ut}^*_{X; Z}$ the image of this map. 
We will frequently identify $\ut_{X; Z}^*$ with $\underline{\ut}_{X; Z}^*$. 

\defi{Let $Z\subset  Z(X)$ be a subset and let $Y\subset \ut_{X; Z}^*$ be a subvariety. 
We say that $Y$ is a {\em set-section} for $\ut_{X; Z}^*$ if $Y$ intersects each $U_X$-orbit of $\ut_{X; Z}^*$ in exactly one point.}
It is clear that a set-section is a section of the natural quotient map $\ut_{X; Z}^*\to \ut_{X; Z}^*/U_X$  but only in set-theoretic terms. 


Next, for two subsets $X, Y$ of vector space we define $X+Y:=\{x+y\mid x\in X, y\in Y\}$. 
In a similar way we define $X-Y$.
For a set $S$ we denote by $|S|$ the number of elements in $S$. 

Let $\ut$ be a finite-dimensional nilpotent Lie algebra with $\vt\subset\ut$ being a subalgebra of codimension~1 (note that $\ut$ is a subfactor of the Lie algebra of upper triangular matrices with the usual restriction on the characteristic $p$). 
We denote by $U, V$ the respective unipotent groups and by $\psi$ the natural map $\ut^*\to\vt^*$. It is clear that $\vt$ is a normal subalgebra of $\ut$ and $V$ is a normal subgroup of $U$. 

Let $S\subset \ut$ be a subset of the center of $\ut$ (this is equivalent to $[S, \ut]=0$). We set $$\ut_{S\ne0}^*:=\{f\in\ut^*\mid \forall s\in S\colon f(s)\ne0\}.$$ 
It is clear that $\ut_{S\ne0}$ is $U$-stable.

The following slight modification of \cite[Proposition 5.11]{IgnatevPetukhov25} holds by the very similar arguments.
\mprop{\label{P:gdid31V0} Let $v_\gamma \in S \setminus(0), v_{\delta}\in \vt, S\subset\vt$. If
$$[v_\gamma, \ut]=0, \hspace{10pt} [v_{\delta}, \vt]=0,\hspace{10pt}[v_{\delta}, \ut]=\Fp_q v_\gamma, \hspace{10pt} [S, \ut]=0~(\text{and~thus~}[S, \vt]=0)$$
then a), b) and c):\\
\textup{a)} $\vt/\Fp_q v_\delta$ is a Lie algebra of dimension $\dim\ut-2=\dim \vt-1$.\\
\textup{b)} Set $S'$ to be the image of $S$ in $\vt/\Fp_q v_{\delta}$ and fix a partial inverse $\phi$ to $\psi$, i.e. a map $\phi: \vt^*\to \ut^*$ such that $$\forall v\in\vt^*\colon \psi(\phi(v))=v.$$ Assume that $Y$ is a set-section for $(\vt/\Fp_q v_\delta)_{S'\ne0}^*$. 
Then $\phi(Y)$ is a set-section for $\ut^*_{S\ne0}$ where we consider $Y$ as a subset of $\vt^*$ via the natural embedding $(\vt/\Fp_q v_\delta)^*\to\vt^*$. \\
\textup{c)} Pick $\lambda\in Y$ and assume that $\mathfrak p\subset \vt$ is a polarization of $\lambda$. Then $\mathfrak p\oplus\Fp_q v_{\delta}$ is a polarization of~$\phi(\lambda)$. This also implies that $\dim (U.\phi(\lambda))=2+\dim (V.\lambda)$.
}{}
We use the previous proposition to prove the next one.
\propp{\label{prop:petukhov}Let $X$ be a quattern and let $\Gamma \subset Z(X)$ be a subset. Suppose that there exist subsets $\Delta, B \subset X \setminus \Gamma$ such that
\begin{enumerate}
\item $|\Delta|=|B|$;
\item $(\Gamma-\Delta)\cap X\subset B$,~where $\Gamma-\Delta=\{\gamma-\delta,\gamma\in\Gamma,\delta\in\Delta\}$;
\item $(B-X)\cap X=\varnothing$,~where $B-X=\{\beta-x,\beta\in B,x\in X\}$;
\item for all $\chi\in \ut_{\Gamma; \Gamma}^*$ the determinant of the matrix $\left( \chi([e_\beta, e_{\delta}])\right)_{\beta\in B, \delta\in\Delta}$ is a nonzero constant.  
\end{enumerate}
Then
\begin{enumerate}
    \item $X\setminus(\Delta\cup B)$ is a C-quattern and $Z(X)\subset Z(X\setminus(\Delta\cup B))$.
    \item Assume that $Y$ is a set-section for $\ut_{X\setminus(\Delta\cup B); \Gamma}^*$. Then $\phi(Y)$ is a set-section for $\ut_{X; \Gamma}^*$\textup, where $\phi$ is the linear map $\ut_{X\setminus(\Delta\cup B)}^*\to\ut_X^*$ defined as follows\textup: $$\phi(e_\tau^*)=e_\tau^*\in\ut_X^*.$$
    \item Pick $\lambda\in Y$ and assume that $\mathfrak p\subset \ut_{X\setminus(\Delta\cup B)}$ is a polarization of $\lambda$. Then $\mathfrak p\oplus\bigoplus_{\delta\in\Delta}\Fp_q e_{\delta}$ is a polarization of~$\phi(\lambda)$. This also implies that $\dim (U_X.\phi(\lambda))=2|\Delta|+\dim (U_{X\setminus(\Delta\cup B)}.\lambda)$.
\end{enumerate}
}{This is essentially  Proposition~\ref{P:gdid31V0} repeated $|\Delta|=|B|$ times. Let $\ut=\ut_X$. 

The first part is trivial. 
We proceed to the second and the third parts. 
Without loss of generality we can assume that the restriction of $Y$ to $\ut_{Z(X)}$ is constant (all our maps can be evaluated separately for the respective pieces). 
We set $\chi$ to be the corresponding linear form $\chi\colon\ut_{\Gamma}\to\Fp_q$.

Set $l:=|\Delta|=|B|$; we also enumerate the elements of $\Delta$ by $\delta_1, \delta_2, \ldots, \delta_l$. 
We set $$\ut_{-i}:=\{x\in\ut\mid \forall i_0\le i\colon \chi([\delta_{i_0}, x])=0\}.$$ 
It is clear that $\ut_{X\setminus (B\cup\Delta)}\subset\ut_{-i}$ for all $i$. 
Thus $\ut_{-i}$ can be completely determined by the matrix $\left( \chi([e_\beta, e_{\delta}])\right)_{\beta\in B, \delta\in\Delta}$. 
Thanks to condition 4) we have $\dim \ut_{-i}=\dim\ut-i$.

Thanks to conditions 1)-4) we see that 
$$\ut\supset \ut_{-1}\supset \ut_{-2}\supset\ldots\supset\ut_{-l} =\ut_{X\setminus B}.$$
Moreover, for any $i$ holds
\begin{enumerate}
    \item $\dim \ut_{-i}=\dim\ut-i$;
    \item $\beta_1, \beta_2,\ldots,\beta_i$ belongs to the center of $\ut_{-i}$;
    \item $\Gamma\subset \ut_{-i}$ for all $i$;
    \item $[\ut_{-i+1}, \beta_i]\subset \ut_{\Gamma}$;
    \item $\chi([\ut_{-i+1}, \beta_i])\ne 0$.    
\end{enumerate}
We set $S:=\{e_\gamma\}_{\gamma\in\Gamma}$. 
Then we apply Proposition~\ref{P:gdid31V0} repeatedly:
$$\ut_{X; \Gamma}^*=\ut^*_{S\ne 0}\to(\ut_{-1}/\Fp_q\beta_1)^*_{S\ne 0}\to (\ut_{-2}/(\Fp_q\beta_1\oplus\Fp_q\beta_2))^*_{S\ne 0}\to$$
$$\to(\ut_{-3}/(\Fp_q\beta_1\oplus\Fp_q\beta_2\oplus\Fp_q\beta_3))^*_{S\ne 0}\to\ldots\to(\ut_{-i}/(\Fp_q\beta_1\oplus\ldots\Fp_q\beta_l))^*_{S\ne 0}\cong \ut_{X\setminus (B\cup\Delta); Z}^*.$$
This finishes the proof.
}

\sect{A classification of subregular orbits}\label{sect:sea_battle_pol}

Before formulating the main theorem of this section, we need to define some subsets of the set of positive roots. More precisely, for $\Phi=B_n$ or $D_n$, put
$$D^{B_n}=\{\epsi_1+\epsi_2,\epsi_3+\epsi_4,\ldots,\epsi_{k-1}+\epsi_k,\epsi_1-\epsi_2,\epsi_3-\epsi_4,\ldots,\epsi_{k-1}-\epsi_k\}\subset\Phi^+,$$
where $k=n-2$ if $n$ is even and $k=n-1$ if $n$ is odd, and
$$D^{D_n}=\{\epsi_1+\epsi_2,\epsi_3+\epsi_4,\ldots,\epsi_{k-1}+\epsi_k,\epsi_1-\epsi_2,\epsi_3-\epsi_4,\ldots,\epsi_{k-1}-\epsi_k\}\subset\Phi^+,$$
where $k=n$ if $n$ is even and $k=n-1$ if $n$ is odd.

Next, in the case, when $\Phi=D_n$ we will denote the set $D^{D_n}$ by $D_\text{reg}$. And in the case $\Phi=B_n$, we will denote by $D_\text{reg}$ one of the sets $D^{B_n}\cup\{\epsi_{n-1}+\epsi_n,\epsi_{n-1}-\epsi_n\}$ or $D^{B_n}\cup\{\epsi_{n-1}\}$ if $n$ is even, and $D^{B_n}\cup\{\epsi_n\}$ if $n$ is odd. Such sets $D_\text{reg}$ of positive roots we will call regular. A set $\xi=(\xi_{\beta})_{\beta\in D_\text{reg}}$ is defined as above, except that $\xi_{\epsi_i-\epsi_{i+1}}$ for $i$ from 1 to $n-1$ and $\xi_{\epsi_n}$ are allowed to be zero.

At the next step we define a set $\Dsr$ by induction on $n$. Firstly, for $B_n$, where $n=2,3$ or 4 we define $\Dsr$ to be one of the following subsets of the positive roots.

If $n=2$, then
$\Dsr=\{\epsi_1-\epsi_2,\epsi_2\}$, where $\xi_{\epsi_1-\epsi_2}$ and $\xi_{\epsi_2}$ can be zero.

If $n=3$, then $\Dsr$ is one of the following subsets:
\begin{itemize}
    \item $\{\epsi_1-\epsi_3,\epsi_2-\epsi_3,\epsi_2,\epsi_1+\epsi_3\}$, where $\xi_{\epsi_1-\epsi_3}$, $\xi_{\epsi_2}$ and $\xi_{\epsi_2-\epsi_3}$ can be zero.
    \item $\{\epsi_2-\epsi_3,\epsi_1,\epsi_2+\epsi_3\}$, where $\xi_{\epsi_2-\epsi_3}$ and $\xi_{\epsi_2+\epsi_3}$ can be zero.
    \item $\{\epsi_2+\epsi_3,\epsi_1-\epsi_3\}$.
\end{itemize}

if $n=4$, then $\Dsr$ is one of the following subsets:
\begin{itemize}
    \item $\{\epsi_1-\epsi_2,\epsi_3-\epsi_4,\epsi_4,\epsi_1+\epsi_2\}$, where $\xi_{\epsi_1-\epsi_2}$, $\xi_{\epsi_3-\epsi_4}$ and $\xi_{\epsi_4}$ can be zero.
    \item $\{\epsi_1-\epsi_3,\epsi_2-\epsi_3,\epsi_2-\epsi_4,\epsi_1+\epsi_3,\epsi_2+\epsi_4\}$, where $\xi_{\epsi_1-\epsi_3}$, $\xi_{\epsi_2-\epsi_3}$ and $\xi_{\epsi_2-\epsi_4}$ can be zero.
    \item $\{\epsi_1-\epsi_3,\epsi_2-\epsi_3,\epsi_1+\epsi_3,\epsi_2\}$, where $\xi_{\epsi_1-\epsi_3}$ and $\xi_{\epsi_2-\epsi_3}$ can be zero.
    \item $\{\epsi_2-\epsi_3,\epsi_1-\epsi_4,\epsi_1+\epsi_4,\epsi_2+\epsi_3\}$, where $\xi_{\epsi_2-\epsi_3}$ can be zero.
    \item $\{\epsi_2-\epsi_3,\epsi_1,\epsi_2+\epsi_3\}$, where $\xi_{\epsi_2-\epsi_3}$ can be zero.
    \item$\{\epsi_2-\epsi_4,\epsi_1+\epsi_4,\epsi_2+\epsi_3\}$.
\end{itemize}

When $n>4$, and $\Phi=B_n$ we define a set $\Dsr$ in one of the following ways. At first, consider a quattern $X=\Phi^+$ and set $$I_1=\{\epsi_1-\epsi_3,\ldots,\epsi_1-\epsi_n,\epsi_1,\epsi_1+\epsi_3,\ldots,\epsi_1+\epsi_n\},$$
$$I_2=\{\epsi_2-\epsi_3,\ldots,\epsi_2-\epsi_n,\epsi_2,\epsi_2+\epsi_3,\ldots,\epsi_2+\epsi_n\}.$$
It can be checked that $Y=X\setminus (I_1\cup I_2)$ is a quattern and it is equal to $\Phi'^+\cup\{\epsi_1-\epsi_2,\epsi_1+\epsi_2\}$, where $\Phi'$ is isomorphic to the root system of type $B_{n-2}$. Now, take an arbitrary $\Dsr\subset\Phi'^+$, $\xi'=(\xi'_{\alpha})_{\alpha\in\Dsr}$ and consider $D=\Dsr\cup\{\epsi_1-\epsi_2,\epsi_1+\epsi_2\}\subset Y$. We also allow the constant $\xi_{\epsi_1-\epsi_2}$ to be zero in that case and put $\xi_{\alpha}=\xi'_{\alpha}$ if $\alpha\in\Dsr$ for $\xi=(\xi_{\beta})_{\beta\in D}$. 

At second, consider the quattern $X=\Phi^+\setminus\{\epsi_1+\epsi_2\}$ and set $$I_1=\{\epsi_1-\epsi_4,\ldots,\epsi_1-\epsi_n,\epsi_1,\epsi_1+\epsi_4,\ldots,\epsi_1+\epsi_n,\epsi_1-\epsi_2\},$$
$$I_2=\{\epsi_3-\epsi_4,\ldots,\epsi_3-\epsi_n,\epsi_3,\epsi_3+\epsi_4,\ldots,\epsi_3+\epsi_n,\epsi_2+\epsi_3\}.$$
It can be checked that $Y=X\setminus (I_1\cup I_2)$ is a quattern and it is equal to $\Phi'^+\cup\{\epsi_1-\epsi_3,\epsi_2-\epsi_3,\epsi_1+\epsi_3\}$, where $\Phi'$ is isomorphic to the root system of type $B_{n-2}$. Now, take an arbitrary $D_\text{reg}\subset\Phi'^+$, $\xi'=(\xi'_{\alpha})_{\alpha\in\Dr}$ and consider $D=\Dr\cup\{\epsi_1-\epsi_3,\epsi_2-\epsi_3,\epsi_1+\epsi_3\}\subset Y$. We also allow the constants $\xi_{\epsi_1-\epsi_3}$ and $\xi_{\epsi_2-\epsi_3}$ to be zero in that case and put $\xi_{\alpha}=\xi'_{\alpha}$ if $\alpha\in\Dr$ for $\xi=(\xi_{\beta})_{\beta\in D}$.

At third, consider a quattern $X=\Phi^+\setminus\{\epsi_1+\epsi_2,\epsi_1+\epsi_3\}$ and set 
$$I_1=\{\epsi_1-\epsi_5,\ldots,\epsi_1-\epsi_n,\epsi_1,\epsi_1+\epsi_5,\ldots,\epsi_1+\epsi_n,\epsi_1-\epsi_3,\epsi_2-\epsi_4,\epsi_2+\epsi_4\},$$
$$I_2=\{\epsi_4-\epsi_5,\ldots,\epsi_4-\epsi_n,\epsi_4,\epsi_4+\epsi_5,\ldots,\epsi_4+\epsi_n,\epsi_1-\epsi_2,\epsi_3-\epsi_4,\epsi_3+\epsi_4\}.$$
It can be checked that $Y=X\setminus (I_1\cup I_2)$ is a quattern and it is equal to $\Phi'^+\cup\{\epsi_1-\epsi_4,\epsi_1+\epsi_4\}$, where $\Phi'$ is isomorphic to the root system of type $B_{n-2}$. Now, take an arbitrary $D_\text{reg}\subset\Phi'^+$, $\xi'=(\xi'_{\alpha})_{\alpha\in\Dr}$ and consider $D=\Dr\cup\{\epsi_1-\epsi_4,\epsi_1+\epsi_4\}\subset Y$. We also put $\xi_{\alpha}=\xi'_{\alpha}$ if $\alpha\in\Dr$ for $\xi=(\xi_{\beta})_{\beta\in D}$.

And finally, consider a subset of the positive roots $X=\Phi^+\setminus\{\epsi_1+\epsi_2,\epsi_1+\epsi_3,\epsi_1-\epsi_4\}$ and set 
$$I_1=\{\epsi_1-\epsi_5,\ldots,\epsi_1-\epsi_n,\epsi_1,\epsi_1+\epsi_5,\ldots,\epsi_1+\epsi_n,\epsi_1-\epsi_3,\epsi_1-\epsi_2,\},$$
$$I_2=\{\epsi_4-\epsi_5,\ldots,\epsi_4-\epsi_n,\epsi_4,\epsi_4+\epsi_5,\ldots,\epsi_4+\epsi_n,\epsi_2+\epsi_4,\epsi_3+\epsi_4\},$$
$$I_3=\{\epsi_2-\epsi_5,\ldots,\epsi_2-\epsi_n,\epsi_2,\epsi_2+\epsi_5,\ldots,\epsi_2+\epsi_n\},$$
$$I_4=\{\epsi_3-\epsi_5,\ldots,\epsi_3-\epsi_n,\epsi_3,\epsi_3+\epsi_5,\ldots,\epsi_3+\epsi_n\}.$$
It can be checked that $Y=X\setminus (I_1\cup I_2\cup I_3\cup I_4)$ is a quattern and it is equal to $\Phi'^+\cup\{\epsi_1+\epsi_4,\epsi_2+\epsi_3\}\cup\Phi''^+$, where $\Phi'$ is isomorphic to the root system of type $B_{n-4}$ and $\Phi''$ is isomorphic to $A_2$. Now, take an arbitrary $\Dr\subset\Phi'^+$, $\Dr'\subset\Phi''^+$, $\xi'=(\xi'_{\alpha})_{\alpha\in\Dr}$, $\xi''=\xi''_{\epsi_2-\epsi_4}$ and consider $D=\Dr\cup\{\epsi_1+\epsi_4,\epsi_2+\epsi_3\}\cup\Dr'\subset Y$ (here, $\Dr'$ is a regular subset for the root system of type $A_2$, in fact, it will be just the root $\epsi_2-\epsi_4$). We also put $\xi_{\alpha}=\xi'_{\alpha}$ if $\alpha\in\Dr$ and $\xi''_{\epsi_2-\epsi_4}=\xi_{\epsi_2-\epsi_4}$ for $\xi=(\xi_{\beta})_{\beta\in D}$.

\defi{We will call a set $D$ constructed in one of the ways described above a subregular set and write $\Dsr$ for the root system $B_n$ and $n>4$.}

We define a set $\Dsr$ in the similar way for the groups of type $D_n$. At first, for $n=3$ or 4, we define a set $\Dsr$ to be one of the following sets.

If $n=3$, then $\Dsr$ is one of the following subsets:
\begin{itemize}
    \item $\{\epsi_1-\epsi_3,\epsi_2-\epsi_3,\epsi_1+\epsi_3\}$, where $\xi_{\epsi_1-\epsi_3}$ and $\xi_{\epsi_2-\epsi_3}$ can be zero.
    \item $\{\epsi_1-\epsi_3,\epsi_2+\epsi_3\}$, where $\xi_{\epsi_2+\epsi_3}$ can be zero.
\end{itemize}

If $n=4$, then $\Dsr$ is one of the following subsets:
\begin{itemize}
    \item $\{\epsi_2-\epsi_3,\epsi_1+\epsi_4,\epsi_1-\epsi_4,\epsi_2+\epsi_3\}$, where $\xi_{\epsi_2-\epsi_3}$ can be zero.
    \item $\{\epsi_1-\epsi_3,\epsi_2-\epsi_3,\epsi_2-\epsi_4,\epsi_1+\epsi_3,\epsi_2+\epsi_4\}$, where $\xi_{\epsi_1-\epsi_3}$, $\xi_{\epsi_2-\epsi_3}$, $\xi_{\epsi_2-\epsi_4}$ and $\xi_{\epsi_2+\epsi_4}$ can be zero.
    \item $\{\epsi_1-\epsi_3,\epsi_1-\epsi_4,\epsi_2+\epsi_3\}$, where $\xi_{\epsi_1-\epsi_4}$ can be zero.
   \item $\{\epsi_1-\epsi_4,\epsi_2-\epsi_4,\epsi_1+\epsi_4\}$, where $\xi_{\epsi_1-\epsi_4}$ can be zero.
    \item $\{\epsi_2-\epsi_4,\epsi_1+\epsi_4,\epsi_2+\epsi_3\}$.
    \item$\{\epsi_1-\epsi_4,\epsi_2+\epsi_4\}$.
\end{itemize}

When $n>4$ we define a set $\Dsr$ in one of the following ways. At first, consider a quattern $X=\Phi^+$ and set $$I_1=\{\epsi_1-\epsi_3,\ldots,\epsi_1-\epsi_n,\epsi_1+\epsi_3,\ldots,\epsi_1+\epsi_n\},$$
$$I_2=\{\epsi_2-\epsi_3,\ldots,\epsi_2-\epsi_n,\epsi_2+\epsi_3,\ldots,\epsi_2+\epsi_n\}.$$
It can be checked that $Y=X\setminus (I_1\cup I_2)$ is a quattern and it is equal to $\Phi'^+\cup\{\epsi_1-\epsi_2,\epsi_1+\epsi_2\}$, where $\Phi'$ is isomorphic to the root system of type $D_{n-2}$. Now, take an arbitrary $\Dsr\subset\Phi'^+$, $\xi'=(\xi'_{\alpha})_{\alpha\in\Dsr}$ and consider $D=\Dsr\cup\{\epsi_1-\epsi_2,\epsi_1+\epsi_2\}\subset Y$. We also allow the constant $\xi_{\epsi_1-\epsi_2}$ to be zero in that case and put $\xi_{\alpha}=\xi'_{\alpha}$ if $\alpha\in\Dsr$ for $\xi=(\xi_{\beta})_{\beta\in D}$. 

At second, consider the quattern $X=\Phi^+\setminus\{\epsi_1+\epsi_2\}$ and set $$I_1=\{\epsi_1-\epsi_4,\ldots,\epsi_1-\epsi_n,\epsi_1+\epsi_4,\ldots,\epsi_1+\epsi_n,\epsi_1-\epsi_2\},$$
$$I_2=\{\epsi_3-\epsi_4,\ldots,\epsi_3-\epsi_n,\epsi_3+\epsi_4,\ldots,\epsi_3+\epsi_n,\epsi_2+\epsi_3\}.$$
It can be checked that $Y=X\setminus (I_1\cup I_2)$ is a quattern and it is equal to $\Phi'^+\cup\{\epsi_1-\epsi_3,\epsi_2-\epsi_3,\epsi_1+\epsi_3\}$, where $\Phi'$ is isomorphic to the root system of type $D_{n-2}$. Now, take an arbitrary $D_\text{reg}\subset\Phi'^+$, $\xi'=(\xi'_{\alpha})_{\alpha\in\Dr}$ and consider $D=\Dr\cup\{\epsi_1-\epsi_3,\epsi_2-\epsi_3,\epsi_1+\epsi_3\}\subset Y$. We also allow the constants $\xi_{\epsi_1-\epsi_3}$, $\xi_{\epsi_2-\epsi_3}$ and $\xi_{\epsi_1+\epsi_3}$ to be zero in that case and put $\xi_{\alpha}=\xi'_{\alpha}$ if $\alpha\in\Dr$ for $\xi=(\xi_{\beta})_{\beta\in D}$.

At third, consider a quattern $X=\Phi^+\setminus\{\epsi_1+\epsi_2,\epsi_1+\epsi_3\}$ and set 
$$I_1=\{\epsi_1-\epsi_5,\ldots,\epsi_1-\epsi_n,\epsi_1+\epsi_5,\ldots,\epsi_1+\epsi_n,\epsi_1-\epsi_3,\epsi_2-\epsi_4,\epsi_2+\epsi_4\},$$
$$I_2=\{\epsi_4-\epsi_5,\ldots,\epsi_4-\epsi_n,\epsi_4+\epsi_5,\ldots,\epsi_4+\epsi_n,\epsi_1-\epsi_2,\epsi_3-\epsi_4,\epsi_3+\epsi_4\}.$$
It can be checked that $Y=X\setminus (I_1\cup I_2)$ is a quattern and it is equal to $\Phi'^+\cup\{\epsi_1-\epsi_4,\epsi_1+\epsi_4\}$, where $\Phi'$ is isomorphic to the root system of type $D_{n-2}$. Now, take an arbitrary $D_\text{reg}\subset\Phi'^+$, $\xi'=(\xi'_{\alpha})_{\alpha\in\Dr}$ and consider $D=\Dr\cup\{\epsi_1-\epsi_4,\epsi_1+\epsi_4\}\subset Y$. We also put $\xi_{\alpha}=\xi'_{\alpha}$ if $\alpha\in\Dr$ for $\xi=(\xi_{\beta})_{\beta\in D}$.

And finally, consider a subset of the positive roots $X=\Phi^+\setminus\{\epsi_1+\epsi_2,\epsi_1+\epsi_3,\epsi_1-\epsi_4\}$ and set 
$$I_1=\{\epsi_1-\epsi_5,\ldots,\epsi_1-\epsi_n,\epsi_1+\epsi_5,\ldots,\epsi_1+\epsi_n,\epsi_1-\epsi_3,\epsi_1-\epsi_2,\},$$
$$I_2=\{\epsi_4-\epsi_5,\ldots,\epsi_4-\epsi_n,\epsi_4+\epsi_5,\ldots,\epsi_4+\epsi_n,\epsi_2+\epsi_4,\epsi_3+\epsi_4\},$$
$$I_3=\{\epsi_2-\epsi_5,\ldots,\epsi_2-\epsi_n,\epsi_2+\epsi_5,\ldots,\epsi_2+\epsi_n\},$$
$$I_4=\{\epsi_3-\epsi_5,\ldots,\epsi_3-\epsi_n,\epsi_3+\epsi_5,\ldots,\epsi_3+\epsi_n\}.$$
It can be checked that $Y=X\setminus (I_1\cup I_2\cup I_3\cup I_4)$ is a quattern and it is equal to $\Phi'^+\cup\{\epsi_1+\epsi_4,\epsi_2+\epsi_3\}\cup\Phi''^+$, where $\Phi'$ is isomorphic to the root system of type $D_{n-4}$ and $\Phi''$ is isomorphic to $A_2$. Now, take an arbitrary $\Dr\subset\Phi'^+$, $\Dr'\subset\Phi''^+$, $\xi'=(\xi'_{\alpha})_{\alpha\in\Dr}$, $\xi''=\xi''_{\epsi_2-\epsi_4}$ and consider $D=\Dr\cup\{\epsi_1+\epsi_4,\epsi_2+\epsi_3\}\cup\Dr'\subset Y$ (here, $\Dr'$ is a regular subset for the root system of type $A_2$, in fact, it will be just the root $\epsi_2-\epsi_4$). We also put $\xi_{\alpha}=\xi'_{\alpha}$ if $\alpha\in\Dr$ and $\xi''_{\epsi_2-\epsi_4}=\xi_{\epsi_2-\epsi_4}$ for $\xi=(\xi_{\beta})_{\beta\in D}$.

\defi{We will call a set $D$ constructed in one of the ways described above a subregular set and write $\Dsr$ for the root system $D_n$ and $n>4$.}

\exam{Here we schematically drew some examples of elements $f_{D,\xi}$ for some $D=D_\text{sreg}$ for $\Phi=D_4$ and $B_4$. The symbol $\otimes$ stands for nonzero constants and the symbol $\square$ for an arbitrary constants $\xi_{\beta}$ for $\beta\in D_\text{sreg}$.
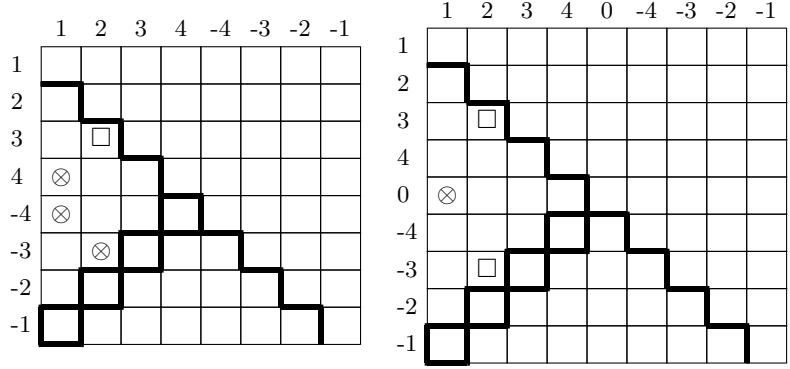
\begin{center}\small
 $\mymatrix{
\lNote{1}\Note{1}\pho& \Note{2}\pho& \Note{3}\pho& \Note{4}\pho& \Note{-4}\pho& \Note{-3}\pho & \Note{-2}\pho & \Note{-1}\\
\lNote{2}\Top{2pt}\Rt{2pt} \pho & \pho& \pho& \pho& \pho& \pho & \pho & \pho\\
\lNote{3}\pho& \Top{2pt}\Rt{2pt}\square& \pho& \pho& \pho& \pho & \pho & \pho\\
\lNote{4}\otimes & \pho& \Top{2pt}\Rt{2pt}\pho& \pho& \pho& \pho & \pho & \pho\\
\lNote{-4}\otimes &\pho & \pho& \Lft{2pt}\Top{2pt}\Bot{2pt}\Rt{2pt} & \pho& \pho& \pho & \pho\\
\lNote{-3}\pho &\otimes& \Lft{2pt}\Top{2pt}\Bot{2pt}\Rt{2pt}\pho& \pho& \Top{2pt}\Rt{2pt}\pho& \pho& \pho & \pho\\
\lNote{-2}\pho &\Lft{2pt}\Top{2pt}\Bot{2pt}\Rt{2pt}& \Top{2pt}\pho& \pho& \pho& \Top{2pt}\Rt{2pt} \pho& \pho & \pho\\
\lNote{-1}\Lft{2pt}\Bot{2pt}\Rt{2pt}\Top{2pt}\pho& \pho& \pho& \pho& \pho& \pho & \Top{2pt}\Rt{2pt}\pho &\pho\\
}
\quad\mymatrix{
\lNote{1}\Note{1}\pho& \Note{2}\pho& \Note{3}\pho& \Note{4}\pho& \Note{0}\pho& \Note{-4}\pho& \Note{-3}\pho & \Note{-2}\pho & \Note{-1}\pho\\
\lNote{2}\Top{2pt}\Rt{2pt}\pho& \pho& \pho& \pho& \pho& \pho & \pho & \pho& \pho\\
\lNote{3}\pho & \square\Top{2pt}\Rt{2pt}& \pho& \pho& \pho& \pho & \pho & \pho& \pho\\
\lNote{4}\pho & \pho& \Top{2pt}\Rt{2pt}\pho& \pho& \pho& \pho & \pho & \pho& \pho\\
\lNote{0}\otimes &\pho & \pho& \Top{2pt}\Bot{2pt}\Rt{2pt}\pho& \pho& \pho& \pho & \pho& \pho\\
\lNote{-4}\pho &\pho & \pho& \Rt{2pt}\Bot{2pt}\Lft{2pt}\pho& \Top{2pt}\Rt{2pt}\pho& \pho& \pho & \pho& \pho\\
\lNote{-3}\pho &\otimes& \Lft{2pt}\Top{2pt}\Bot{2pt}\Rt{2pt} & \pho& \pho& \Top{2pt}\Rt{2pt}\pho& \pho & \pho& \pho\\
\lNote{-2}\pho &\Lft{2pt}\Top{2pt}\Bot{2pt}\Rt{2pt} & \pho& \pho& \pho& \pho& \Top{2pt}\Rt{2pt}\pho & \pho& \pho\\
\lNote{-1}\Lft{2pt}\Top{2pt}\Bot{2pt}\Rt{2pt}& \pho& \pho& \pho& \pho& \pho & \pho &\Top{2pt}\Rt{2pt}\pho& \pho\\
}$
\end{center}
At the left picture we have $D=\{\epsi_1-\epsi_4,\epsi_2+\epsi_3,\epsi_2-\epsi_3,\epsi_1+\epsi_4\}$ and $D=\{\epsi_2-\epsi_3,\epsi_2+\epsi_3,\epsi_1\}$ at the right picture.}
\theop{Let\label{theo_orb_submax_dim} $D=D_\text{sreg}\subset\Phi^+$ be a subregular subset of the positive roots and $\xi$ be as above\textup{,} then put 
\begin{equation*}
f=f_{D, \xi} = \sum_{\beta\in D}\xi_{\beta}e_{\beta}^*\in\ut^*.
\end{equation*}
Then the orbit of $f$ has pre-maximal dimension. Conversely\textup{,} each orbit of pre-maximal dimension contains exactly one such linear form.}{In both cases, $B_n$ and $D_n$, we will proceed by induction on $n$. The base of induction for $B_2$, $B_3$, $B_4$ and $D_3$, $D_4$ can be easily checked. Indeed, it is easy to see that all the orbits described above are subregular, i.e., have the pre-maximal dimension, and they are distinct. On the other hand, the number of such orbits is known, see \cite{GoodwinMoschRoehrle16}. Hence, we can conclude that we described all of them, so we have checked the base of induction. 

Now, assume that $n>4$ and $\Phi=B_n$ or $D_n$. Firstly we want to prove that, for a given linear form $$f=\sum_{\beta\in\Phi^+}\xi_{\beta}e_{\beta}^*\in\ut^*$$ and $\delta\prec\epsi_1+\epsi_4$, where $\delta$ is the biggest root with nonzero constant $\xi_{\delta}$ with respect to the order~$\prec$, the dimension of the orbit of $f$ is less than pre-maximal. Let us assume that $\delta=\epsi_1+\epsi_5$, i.e., $\xi_{\epsi_1+\epsi_2}=\xi_{\epsi_1+\epsi_3}=\xi_{\epsi_1+\epsi_4}=0$. Consider the quattern $Y=\Phi^+\setminus\{\epsi_1+\epsi_2,\epsi_1+\epsi_3,\epsi_1+\epsi_4\}$. Put $Z=\langle e_{\delta}\rangle\subset Z(Y)$ and define a linear form 
\begin{equation}
\wt f=\sum_{\beta\in Y}\xi_{\beta}e_{\beta}^*\in\ut^*_{Y;Z}.
\end{equation}\label{form:projection}

The dimension of the orbit of the form $f$ can be calculated as $\text{rank}(\beta_{f})$ (recall that $\beta_f$ is a bilinear form on $\ut$ defined by $\beta_f(x,y)=f([x,y])$). We want to compare the ranks of the bilinear forms $\beta_f$ and $\beta_{\wt f}$ defined on $\ut^*$ and $\ut^*_{Y}$ respectively. Let us denote the matrix of the form $\beta_f$ by $A$ and the matrix of the form $\beta_{\wt f}$ by $\wt A$ and denote the ideal of $\ut$ spanned by the vectors $e_{\epsi_1+\epsi_2},e_{\epsi_1+\epsi_3},e_{\epsi_1+\epsi_4}$ by $I$ (columns and rows are numbered by the roots from $\Phi^+$ which are sorted by the order $\prec$), so $A_{\alpha,\beta}=\beta_f(e_{\alpha},e_{\beta})\in\Fp_q$.
If $[e_{\alpha},e_{\beta}]\notin I$ for some $\alpha,\beta\in\Phi^+$, then $\beta_f(e_{\alpha},e_{\beta})=\beta_{\wt f}(e'_{\alpha},e'_{\beta})$ (here $e'_{\alpha}$ is the class of the vector $e_{\alpha}$), if $[e_{\alpha},e_{\beta}]\in I$, then $\beta_f(e_{\alpha},e_{\beta})=\beta_{\wt f}(e'_{\alpha},e'_{\beta})=0$ because $f|_{I}=0$ and $[e'_{\alpha},e'_{\beta}]=0$ and finally, if $e_{\alpha}$ or $e_{\beta}$ belongs to $I$, we have $\beta_f(e_{\alpha},e_{\beta})=0$ and the corresponding row and column in $\wt A$ does not even exist. Thus, the matrix $\wt A$ is obtained from $A$ by deleting some zero rows and columns, so their ranks are equal. 


Hence, we can consider the form $\wt f$ and its orbit in $\ut^*_{Y;Z}$ instead of $f$ and its orbit in $\ut^*$. At first, consider the case when $\Phi=B_n$.
Let $Y_1$ be defined as the set $Y\setminus (P_1\cup Q_1)$ for $$P_1=\{\epsi_1-\epsi_6,\epsi_1-\epsi_7,\ldots,\epsi_1-\epsi_n,\epsi_1,\epsi_1+\epsi_6,\epsi_1+\epsi_7,\ldots,\epsi_1+\epsi_n\},$$
$$Q_1=\{\epsi_5-\epsi_6,\epsi_5-\epsi_7,\ldots,\epsi_5-\epsi_n,\epsi_5,\epsi_5+\epsi_6,\epsi_5+\epsi_7,\ldots,\epsi_5+\epsi_n\},$$
obtained from the set $Y$ by step-by-step removing of the roots $\epsi_1+\epsi_i$ and $\epsi_5-\epsi_i$ for $i$ from 6 to $n$ pairwise, then $\epsi_1$ and $\epsi_5$ and finally $\epsi_1-\epsi_i$ and $\epsi_5+\epsi_i$ for $i$ from $n$ down to 6 pairwise (in case, when $n=5$ we assume that the sets $P_1$ and $Q_1$ consist only of $\epsi_1$ and $\epsi_5$, respectively).

Consider the map $\vfi:\ut^*_{Y}\to\ut^*=\ut^*_X$ defined as follows for any quattern $Y$: for $e_{\tau}^*\in\ut_Y^*$,
$$\vfi(e^*_{\tau})=e^*_{\tau}\in\ut^*.$$ According to \cite[Proposition 5.11]{IgnatevPetukhov25},
the set $Y_1$ is a quattern and there is a form $f_1\in\ut^*_{Y_1;Z}$ such that $\vfi(f_1)$ lies in the orbit of the form $\wt f$ in $\ut^*_{Y;Z}$ and $$\dim\Omega_{\wt f}=\dim\Omega_{f_1}+|P_1|+|Q_1|=\dim\Omega_{f_1}+2(2n-9).$$ Now, we want to examine the maximal possible dimension of orbit in $\ut^*_{Y_1;Z}$. Let us assume that $f_1\in\ut^*_{Y_1;Z_1}$, where $Z_1=Z\cup\langle e_{\epsi_2+\e_3}\rangle\cup\langle e_{\e_1-\e_5}\rangle$. Put $Y_2=Y_1\setminus (P_2\cup Q_2)$ for $$P_2=\{\e_2+\e_4,\e_3+\e_4,\e_4-\e_5\},$$ $$Q_2=\{\e_3-\e_4,\e_2-\e_4,\e_1-\e_4\}.$$ Here we remove the roots from the set $Y_1$ pairwise in the strict order from the left to the right (i.e., $\e_2+\e_4$ with $\e_3-\e_4$ etc.). According to \cite[Proposition 5.11]{IgnatevPetukhov25}, the set $Y_2$ is a quattern, and there exists $f_2\in\ut^*_{Y_2;Z_1}$ such that $\vfi(f_2)$ lies in the orbit of the form $f_1$ and  
$$\dim\Omega_{f_1}=\dim\Omega_{f_2}+|P_2|+|Q_2|=\dim\Omega_{f_2}+6.$$

Now, we can apply the Proposition \ref{prop:petukhov}. Precisely, consider the set of roots $$Y_3=Y_2\setminus\{e_2+\e_5,\e_3+\e_5,\e_2-\e_5,e_3-\e_5,\e_1-\e_2,\e_1-\e_3\}.$$ This proposition gives us that $Y_3$ is a quattern and there exists $f_3\in\ut^*_{Y_3;Z_1}$ such that $\vfi(f_3)$ lies in the orbit of the form $f_2$ and  
$\dim\Omega_{f_2}=\dim\Omega_{f_3}+6.$

Finally, we can apply \cite[Proposition 5.11]{IgnatevPetukhov25} to the quattern $Y_3$ and the sets $P_4$ and $Q_4$ of positive roots for $$P_4=\{\epsi_2-\epsi_6,\epsi_2-\epsi_7,\ldots,\epsi_2-\epsi_n,\epsi_2,\epsi_2+\epsi_6,\epsi_2+\epsi_7,\ldots,\epsi_2+\epsi_n\},$$
$$Q_4=\{\epsi_3-\epsi_6,\epsi_3-\epsi_7,\ldots,\epsi_3-\epsi_n,\epsi_3,\epsi_3+\epsi_6,\epsi_3+\epsi_7,\ldots,\epsi_3+\epsi_n\}.$$
More precisely, we put $Y_4=Y_3\setminus(P_4\cup Q_4)$, where $Y_4$ obtained from the set $Y_3$ by step-by-step removing of the roots $\epsi_2+\epsi_i$ and $\epsi_3-\epsi_i$ for $i$ from 6 to $n$ pairwise, then $\epsi_2$ and $\epsi_3$ and finally $\epsi_2-\epsi_i$ and $\epsi_3+\epsi_i$ for $i$ from $n$ down to 6 pairwise (in case, when $n=5$ we assume that the sets $P_4$ and $Q_4$ consist only of $\epsi_2$ and $\epsi_3$, respectively). The set $Y_4$ is a quattern and there exists $f_4\in\ut^*_{Y_4;Z_1}$ such that $\vfi(f_4)$ lies in the orbit of the form $f_3$ and  
$$\dim\Omega_{f_3}=\dim\Omega_{f_4}+|P_4|+|Q_4|=\dim\Omega_{f_4}+2(2n-9).$$

It is easy to check that the quotient algebra $\ut_{Y_4}$ is isomorphic to the sum of its ideals $\ut(B_{n-4})\oplus\langle e_{\e_1+\e_5}\rangle\oplus\langle e_{\epsi_1-\epsi_5}\rangle\oplus\langle e_{\epsi_2-\epsi_3}\rangle\oplus\langle e_{\epsi_4+\epsi_5}\rangle\oplus\langle e_{\epsi_2+\epsi_3}\rangle$. It follows that $\dim\Omega_{f_4}$ in $\ut^*_{Y_4;Z_1}$ is no more than $(n-4)(n-5)$, the maximal possible dimension of an orbit in $\ut^*(B_{n-4})$. So, $\dim\Omega_{f_1}$ in $\ut^*_{Y_1;Z_1}$ is no more than $$(n-4)(n-5)+2(2n-9)+12.$$ Since $\ut^*_{Y_1;Z_1}$ is an open subspace of $\ut^*_{Y_1;Z}$ as well as the set of orbits of maximal dimension, its intersection with the set of orbits of maximal dimension is non-empty, so the maximal dimension of an orbit in $\ut^*_{Y_1,Z}$ is also no more than $(n-4)(n-5)+2(2n-9)+12$. (In this place we can consider the orbit of our linear form over the closure of the field $\Fp_q$ and then take its set of the $\Fp_q$-points. The well-known fact tells us that the obtained set will be the orbit of our linear form over the field $\Fp_q$.) Thus, the maximal possible dimension of an orbit in $\ut^*_{Y;Z}$ is no more than $(n-4)(n-5)+4(2n-9)+12$. Arguing as above, we obtain that the maximal possible dimension of an orbit in $\ut^*_Y$ is also no more than $$(n-4)(n-5)+4(2n-9)+12=n(n-1)-4$$ which is less than pre-maximal dimension of an orbit in $\ut^*=\ut^*(B_n)$. So, if the root $\delta$ is bigger than $\epsi_1+\epsi_4$, we can consider the orbit of $f$ as orbit in the quotient algebra $\ut^*_{Y}$. It concludes the proof of the first part.



Now, we know that a linear form $f=(\xi_{\beta})_{\beta\in\Phi^+}$ having the orbit of maximal dimension has to have nonzero constant $\xi_\delta$, where $\delta$ is the biggest root such that $\xi_{\delta}\ne 0$ with $\delta=\epsi_1+\epsi_2$, $\epsi_1+\epsi_3$ or $\epsi_1+\epsi_4$. We aim to present a set-section of the set of orbits of pre-maximal dimension in each case. Let us consider the case when $\delta=\epsi_1+\epsi_2$ at first. We proceed by induction on $n$. Firstly, put $X=\Phi^+$, $\delta=\epsi_1+\epsi_2$, $Z=\langle\delta\rangle= Z(X)$, and $$I_1=\{\epsi_1-\epsi_3,\ldots,\epsi_1-\epsi_n,\epsi_1,\epsi_1+\epsi_3,\ldots,\epsi_1+\epsi_n\},$$
$$I_2=\{\epsi_2-\epsi_3,\ldots,\epsi_2-\epsi_n,\epsi_2,\epsi_2+\epsi_3,\ldots,\epsi_2+\epsi_n\}.$$ We define the set $Y=X\setminus(I_1\cup I_2)$ to be the set obtained by step-by-step removing of the roots $\epsi_1+\epsi_i$ and $\epsi_2-\epsi_i$ for $i$ from 3 to $n$ pairwise, then $\epsi_1$ and $\epsi_2$ and finally $\epsi_1-\epsi_i$ and $\epsi_2+\epsi_i$ for $i$ from $n$ down to 3 pairwise.
According to \cite[Proposition 5.11]{IgnatevPetukhov25},
the set $Y$ is a quattern and there is a form $f'\in\ut^*_{Y;Z}$ such that $\vfi(f')$ lies in the orbit of the form $\wt f$ in $\ut^*_{Y;Z}$ for $$\wt f=\sum_{\beta\in X}\xi_{\beta}e_{\beta}^*\in\ut^*_{X;Z},$$ and $$\dim\Omega_{\wt f}=\dim\Omega_{f'}+|I_1|+|I_2|=\dim\Omega_{f'}+2(2n-3).$$
It is easy to check that the quotient algebra $\ut_{Y}$ is isomorphic to the sum of its ideals $\ut(B_{n-2})\oplus\langle e_{\e_1+\e_2}\rangle\oplus\langle e_{\epsi_1-\epsi_2}\rangle$. Now, if we assume that the orbit of $\wt f$ has the pre-maximal dimension $n(n-1)-2$, we obtain that the orbit of $f'$ has to have the dimension $$n(n-1)-2-2(2n-3)=(n-2)(n-3)-2,$$ i.e., the pre-maximal dimension of an orbit in $\ut^*(B_{n-2})$. This concludes the proof in this case.  

Now, we move to the case, when $\delta=\e_1+\e_3$. Firstly, put $X=\Phi^+\setminus\{\e_1+\e_2\}$, $Z=\langle\delta\rangle= Z(X)$, and $$I_1=\{\epsi_1-\epsi_4,\ldots,\epsi_1-\epsi_n,\epsi_1,\epsi_1+\epsi_4,\ldots,\epsi_1+\epsi_n,\e_2+\e_3\},$$
$$I_2=\{\epsi_3-\epsi_4,\ldots,\epsi_3-\epsi_n,\epsi_3,\epsi_3+\epsi_4,\ldots,\epsi_3+\epsi_n,\e_1-\e_2\}.$$ We define the set $Y=X\setminus(I_1\cup I_2)$ to be the set obtained by step-by-step removing of the roots $\epsi_1+\epsi_i$ and $\epsi_3-\epsi_i$ for $i$ from 4 to $n$ pairwise, then $\epsi_1$ and $\epsi_3$, then $\epsi_1-\epsi_i$ and $\epsi_3+\epsi_i$ for $i$ from $n$ down to 4, and finally $\e_2+\e_3$ and $\e_1-\e_2$ pairwise.
According to \cite[Proposition 5.11]{IgnatevPetukhov25},
the set $Y$ is a quattern and there is a form $f'\in\ut^*_{Y;Z}$ such that $\vfi(f')$ lies in the orbit of the form $\wt f$ in $\ut^*_{Y;Z}$ for $$\wt f=\sum_{\beta\in X}\xi_{\beta}e_{\beta}^*\in\ut^*_{X;Z},$$ and $$\dim\Omega_{\wt f}=\dim\Omega_{f'}+|I_1|+|I_2|=\dim\Omega_{f'}+2(2n-4).$$
It is easy to check that the quotient algebra $\ut_{Y}$ is isomorphic to the sum of its ideals $\ut(B_{n-2})\oplus\langle e_{\e_1+\e_3}\rangle\oplus\langle e_{\epsi_1-\epsi_3}\rangle\oplus\langle e_{\epsi_2+\epsi_3}\rangle$. Now, if we assume that the orbit of $\wt f$ has the pre-maximal dimension $n(n-1)-2$, we obtain that the orbit of $f'$ has to have the dimension $$n(n-1)-2-2(2n-4)=(n-2)(n-3),$$ i.e., the maximal dimension of an orbit in $\ut^*(B_{n-2})$. This concludes the proof in this case. 

Now, we move to the case, when $\delta=\e_1+\e_4$. Firstly, put $X=\Phi^+\setminus\{\e_1+\e_2,\e_1+\e_3\}$, $Z=\langle\delta\rangle=Z(X)$, and $$I_1=\{\epsi_1-\epsi_5,\ldots,\epsi_1-\epsi_n,\epsi_1,\epsi_1+\epsi_5,\ldots,\epsi_1+\epsi_n\},$$
$$I_2=\{\epsi_4-\epsi_5,\ldots,\epsi_4-\epsi_n,\epsi_4,\epsi_4+\epsi_5,\ldots,\epsi_4+\epsi_n\}.$$ We define the set $Y=X\setminus(I_1\cup I_2)$ to be the set obtained by step-by-step removing of the roots $\epsi_1+\epsi_i$ and $\epsi_4-\epsi_i$ for $i$ from 5 to $n$ pairwise, then $\epsi_1$ and $\epsi_4$, then $\epsi_1-\epsi_i$ and $\epsi_4+\epsi_i$ for $i$ from $n$ down to 5 pairwise. 
According to \cite[Proposition 5.11]{IgnatevPetukhov25},
the set $Y$ is a quattern and there is a form $f'\in\ut^*_{Y;Z}$ such that $\vfi(f')$ lies in the orbit of the form $\wt f$ in $\ut^*_{Y;Z}$ for $$\wt f=\sum_{\beta\in X}\xi_{\beta}e_{\beta}^*\in\ut^*_{X;Z},$$ and $$\dim\Omega_{\wt f}=\dim\Omega_{f'}+|I_1|+|I_2|=\dim\Omega_{f'}+2(2n-7).$$

Here we have some separate cases. At first, we can assume that $\xi_{\e_1-\e_4}\neq0$ and $\xi_{\e_2+\e_3}\neq0$. Put 
\begin{equation*}
J_1=\{\epsi_1-\epsi_3,\epsi_2-\epsi_4,\epsi_2+\epsi_4\}\text{, and } 
J_2=\{\epsi_1-\epsi_2,\epsi_3-\epsi_4,\epsi_3+\epsi_4\}.
\end{equation*}
Now, we can apply Proposition \ref{prop:petukhov}. Precisely, for $Y'=Y\setminus(J_1\cup J_2)$ we have that $Y'$ is a quattern and there is a linear form $f''\in\ut^*_{Y',Z'}$ for $Z'=Z\cup\{\e_1-\e_4,\e_2+\e_3\}$ such that $\vfi(f'')$ lies in the orbit of $f'$ and $$\dim\Omega_{f'}=\dim\Omega_{f''}+6.$$

It is easy to check that the quotient algebra $\ut_{Y'}$ is isomorphic to the sum of its ideals $\ut(B_{n-2})\oplus\langle e_{\e_1+\e_4}\rangle\oplus\langle e_{\epsi_1-\epsi_4}\rangle\oplus\langle e_{\epsi_2+\epsi_3}\rangle$. Now, if we assume that the orbit of $\wt f$ has the pre-maximal dimension $n(n-1)-2$, we obtain that the orbit of $f''$ has to have the dimension $$n(n-1)-2-2(2n-7)-6=(n-2)(n-3),$$ i.e., the maximal dimension of an orbit in $\ut^*(B_{n-2})$. This concludes the proof in this case.

At second, assume that $\xi_{\e_1-\e_4}=0$ and $\xi_{\e_2+\e_3}\neq0$. Put $$J_1=\{\epsi_2-\epsi_5,\ldots,\epsi_2-\epsi_n,\epsi_2,\epsi_2+\epsi_5,\ldots,\epsi_2+\epsi_n,\e_1-\e_3,\e_1-\e_2\},$$
$$J_2=\{\epsi_3-\epsi_5,\ldots,\epsi_3-\epsi_n,\epsi_3,\epsi_3+\epsi_5,\ldots,\epsi_3+\epsi_n,\e_3+\e_4,\e_2+\e_4\}.$$
We define the set $Y'=Y\setminus(J_1\cup J_2)$ to be the set obtained by step-by-step removing of the roots $\epsi_2+\epsi_i$ and $\epsi_3-\epsi_i$ for $i$ from 5 to $n$ pairwise, then $\epsi_2$ and $\epsi_3$, then $\epsi_2-\epsi_i$ and $\epsi_3+\epsi_i$ for $i$ from $n$ down to 5 and finally $\e_1-\e_3$ with $\e_3+\e_4$ and $\e_1-\e_2$ with $\e_2+\e_4$ pairwise. 
According to \cite[Proposition 5.11]{IgnatevPetukhov25},
the set $Y'$ is a quattern and there is a form $f''\in\ut^*_{Y';Z'}$ for $Z'=Z\cup\{\e_2+\e_3\}$ such that $\vfi(f'')$ lies in the orbit of $f'$ and $$\dim\Omega_{ f'}=\dim\Omega_{f''}+|J_1|+|J_2|=\dim\Omega_{f''}+2(2n-7)+4.$$

It is easy to check that the quotient algebra $\ut_{Y'}$ is isomorphic to the sum of its ideals $\ut(B_{n-4})\oplus\langle e_{\e_1+\e_4}\rangle\oplus\langle e_{\epsi_2+\epsi_3}\rangle\oplus\ut(A_2)$. Now, if we assume that the orbit of $\wt f$ has the pre-maximal dimension $n(n-1)-2$, we obtain that the orbit of $f''$ has to have the dimension $$n(n-1)-2-4(2n-7)-4=(n-4)(n-5)+2,$$ i.e., the maximal dimension of an orbit in $\ut^*(B_{n-4})\oplus\ut^*(A_2)$. This concludes the proof in this case. Arguing in the same way, we can easily obtain that the form $\wt f$ can not have the pre-maximal dimension in the case when $\xi_{\e_2+\e_3}=0$.

A proof for the case when $\Phi=D_n$ can be conducted in the exact same way as for $B_n$ with only difference that roots $\e_i$ for $i$ from 1 to $n$ do not exist, so we will not add them to the sets $P$, $Q$, $I$ and~$J$.
}

As a corollary, we can prove Isaacs' conjecture in that case.
\corop{The number\label{coro:Isaacs_C} of subregular orbits of group of type $B_n$ and $D_n$ over the finite field $\Fp_q$ is a polynomial in $v=q-1$ with integer non-negative coefficients.}{For $n>4$, from the proof of the theorem above we see that such a number $S_n(\Phi)$ is given by the recurrent condition $$S_n(\Phi)=v(v+1)S_{n-2}(\Phi)+(v^3+3v^2+v)S'_{n-2}(\Phi)+v^3S'_{n-4}(\Phi),$$ where $S'_{n}(\Phi)$ is the expression for the number of regular orbits of group of type $B_n$ or $D_n$ (see \cite{IgnatevVenchakov26}, for example). It is obvious from the construction of such orbits that the base of induction for $n
\leq 4$ gives us that $S_2(\Phi)$, $S_3(\Phi)$ and $S_4(\Phi)$ for $\Phi=B_n$ or $D_n$ are polynomials in $v$. So, assuming that $S_{n-2}(\Phi)$, $S'_{n-2}(\Phi)$ and $S'_{n-4}(\Phi)$ are polynomials in $v$ with integer non-negative coefficients, we obtain the required statement.}

\bigskip\textsc{Mikhail Ignatev: National Research University Higher School of Economics,\break Po\-krov\-sky Boulevard 11, 109028, Moscow, Russia}

\emph{E-mail address}: \texttt{mihail.ignatev@gmail.com}

\bigskip\textsc{Mikhail Venchakov: National Research University Higher School of Economics, Pokrovsky Boulevard 11, 109028, Moscow, Russia}

\emph{E-mail address}: \texttt{mihail.venchakov@gmail.com}


\begin{thebibliography}{XXXXX}
\bibitem[Bo02]{Bourbaki03} N. Bourbaki. Lie groups and Lie algebras. Chapters 4--6. Elements of mathematics (Berlin). Springer--Verlag, Berlin, 2002.
\bibitem[An95i]{Andre95i} C.A.M. Andre. Basic characters of the
unitriangular group. J. Algebra \textbf{175} (1995), 287--319.
\bibitem[An95ii]{Andre95ii} C.A.M. Andre. Basic sums of coadjoint orbits of the
unitriangular group. J. Algebra \textbf{176} (1995), 959--1000.
\bibitem[AN06]{AndreNeto06} C.A.M. Andre, A.M. Neto. Super-characters of finite
unipotent groups of types $B_n$, $C_n$ and $D_n$. J. Algebra \textbf{305} (2006) 394--429.

\bibitem[GMR16]{GoodwinMoschRoehrle16} S.M. Goodwin, P. Mosch, G. R\"ohrle. On the coadjoint orbits of maximal unipotent subgroups of reductive groups. Transformation Groups \textbf{21} (2016), 399--426.

\bibitem[Hi60]{Higman60} G. Higman. Enumerating $p$-groups. I. Inequalities. Proc. London Math. Soc. (3) \textbf{10} (1960), 24--30.


\bibitem[IP09]{IgnatevPanov09} M.V. Ignatev, A.N. Panov. Coadjoint orbits of the group $\mathrm{UT}(7, K)$. J. Math. Sci. \textbf{156} (2009), no. 2, 292--312.

\bibitem[IP25]{IgnatevPetukhov25} M. Ignatev, A. Petukhov. Coadjoint orbits of low dimension for nilradicals of Borel subalgebras in classical types. Expositiones Mathematicae, accepted, arXiv: \texttt{math.RT/2507.20332}.

\bibitem[IV26]{IgnatevVenchakov26} M. Ignatev, M. Venchakov. Orbits of maximal and submaximal dimension for Sylow $p$-subgroups of finite classical groups. Mat. Sbornik, submitted, arXiv: \texttt{arXiv:2608.27523} 

\bibitem[Is07]{Isaacs07} I.M. Isaacs. Counting characters of upper triangular groups. J. Algebra \textbf{315} (2007),\break 698--719.

\bibitem[Ka77]{Kazhdan77} D. Kazhdan. Proof of Springer's hypothesis.
Israel J. Math. \textbf{28} (1977), 272--286.

\bibitem[Ki62]{Kirillov62} A.A. Kirillov. Unitary representations of nilpotent Lie groups. Russian Math. Surveys \textbf{17}~(1962), 53--104.

\bibitem[Ki04]{Kirillov04} A.A. Kirillov. Lectures on the orbit method. Grad. Stud. in Math. \textbf{64}, AMS, 2004.


\bibitem[Ko12]{Kostant12} B. Kostant. The cascade of orthogonal roots and the coadjoint structure of the nilradical of a Borel subgroup of a semisimple Lie group. Moscow Math. J. \textbf{12} (2012), no. 3, 605--620.

\bibitem[Ko13]{Kostant13} B. Kostant. Center of $U(\nt)$, cascade of orthogonal roots and a construction of Lipsman--Wolf. In: Lie groups: structure, actions and representations, Progr. in Math. \textbf{306}. Birkh$\ddot{\mathrm{a}}$user, 2013, 163--174.


\bibitem[Le10]{Le10} T. Le. Counting irreducible representations of large degree of the upper triangular groups. J. Algebra \textbf{324} (2010), 1803--1817.

\bibitem[Le74]{Lehrer74} G.I. Lehrer. Discrete series and the unipotent subgroup. Compositio
Math. \textbf{28} (1974),\break fasc.~1, 9--19.

\bibitem[Lou11]{Loukaki11} M. Loukaki. Counting characters of small degree in upper unitriangular groups. J. Pure Appl. Algebra \textbf{215} (2011), no. 2, 154--160.


\bibitem[Mb11]{Marberg11} E. Marberg. Combinatorial methods of character enumeration for the unitriangular group. J. Algebra \textbf{345} (2011), 295--323.


\bibitem[Mj97]{Marjoram97} M. Marjoram. Irreducible characters of Sylow $p$-subgroups of the orthogonal group. Extracta Mathematicae \textbf{12} (1997), no. 3, 315--319.

\bibitem[Mj97']{Marjoram97'} M. Marjoram. Irreducible characters of Sylow p-subgroups of classical groups. PhD thesis,
National University of Ireland, Dublin, 1997.

\bibitem[Mj99]{Marjoram99} M. Marjoram. Irreducible characters of small degree of the unitriangular group. Irish Math. Soc. Bull. \textbf{42} (1999), 21--31.




\bibitem[Ve26]{Venchakov} M. Venchakov. Rook placements and coadjoint orbits for maximal unipotent subgroups of finite symplectic groups, arXiv: \texttt{math.RT/2602.14933}.

\bibitem[VLA03]{VeraLopezArregi03} A. Vera-L\'opez, J.M. Arregi. Conjugacy classes in unitriangular matrices. Linear Algebra and its Applications \textbf{370} (2003), 85--124.



\end{thebibliography}

\begin{thebibliography}{99}

\bibitem{Humpreys1}Хамфри Дж. Введение в теорию алгебр Ли and их
представлений. --- М.: МЦНМО, 2003.
\bibitem{Ignatev}Игнатьев М.В. Субрегулярные характеры
унитреугольной группы над конечным полем. Фундаментальная and
прикладная математика, т. \textbf{13}, вып. \textbf{5}, 2007, с.
103-125, см. также arXiv:math.RT/0801.3079v2.
\bibitem{IgnatevPanov} Игнатьев М.В., Панов А.Н.
Коприсоединённые орбиты группы $\mathrm{UT}(7, K)$. Фундаментальная
and прикладная математика, т. \textbf{13}, вып. \textbf{5}, 2007, с.
127-159, см. также arXiv: math.RT/0603649v3.
\bibitem{Kirillov1} Кириллов А.А. Лекции по методу орбит. ---
Новосибирск: Научная книга ИДМИ, 2002.
\bibitem{Kirillov2} Кириллов А.А. Унитарные представления
нильпотентных групп Ли. УМН, т. \textbf{17}, 1962, с. 57-110.
\bibitem{Kirillov3} Кириллов А.А. Метод орбит and конечные группы. --- М.: МЦНМО, МК НМУ, 1998.
\bibitem{Panov} Панов А.Н. Инволюции в $S_n$ and ассоциированные
коприсоединённые орбиты. Зап. научн. сем. ПОМИ, т. \textbf{349},
2007, с. 150-173, см. также arXiv: math.RT/0801.3022v1.
\end{thebibliography}
\end{document}